\documentclass[leqno,12pt]{amsart}

\usepackage{amssymb}
\usepackage{amscd}              
\usepackage[all]{xy}            
\xyoption{curve}
\CompileMatrices

\newcommand{\Aa}{{\mathbb{A}}}

\newcommand{\bS}{\bar{S}}

\newcommand{\C}{\mathbb{C}}

\newcommand{\D}{\mathcal{D}}
\newcommand{\diag}{\operatorname{diag}}        
\newcommand{\Dt}{\tilde{D}}
\newcommand{\dual}{\lor}

\newcommand{\GL}{\operatorname{GL}}

\newcommand{\Gm}{\mathbb{G}_m}

\newcommand{\injto}{\hookrightarrow}

\newcommand{\isomorph}{\cong}           
\newcommand{\isomto}{\overset{\sim}{\to}}

\newcommand{\KGL}{\operatorname{KGL}}

\newcommand{\Lie}{\operatorname{Lie}}
\newcommand{\m}{{\mathfrak{m}}}                  

\newcommand{\Mt}{\tilde{M}}

\newcommand{\Oo}{{\mathcal{O}}}                  

\newcommand{\Q}{{\mathbb{Q}}}                    

\newcommand\rk{\operatorname{rk}}

\newcommand{\Spec}{\operatorname{Spec}\, }       

\newcommand{\Sym}{\operatorname{Sym}\, }
\newcommand{\tensor}{\otimes}

\newcommand{\V}{{\mathcal V}}

\newcommand{\X}{{\mathcal{X}}}
\newcommand{\Xt}{\tilde{X}}

\newtheorem{theorem}{Theorem}[section]
\newtheorem{proposition}[theorem]{Proposition}

\newtheorem{corollary}[theorem]{Corollary}

\theoremstyle{definition}
\newtheorem{definition}[theorem]{Definition}

\newtheorem{example}[theorem]{Example}

\begin{document}

\title[Vanishing of top equivariant Chern classes of regular embeddings]
{Vanishing of top equivariant Chern classes \\ of regular embeddings}
\author[M. Brion and I. Kausz]{M. Brion and I. Kausz}
\thanks{The second author has been partially supported by the DFG} 
\date{August 25, 2005}

\address{Michel Brion, Universit\'e de Grenoble I, Institut Fourier,
BP 74, 38402 Saint-Martin d'H\`eres, France}
\email{Michel.Brion@ujf-grenoble.fr}

\address{Ivan Kausz, NWF I - Mathematik, Universit\"{a}t Regensburg, 
93040 Regensburg, 
Germany}
\email{ivan.kausz@mathematik.uni-regensburg.de}

\maketitle



\section{Introduction}

Let $G$ be a connected affine algebraic group and let
$X$ be a regular $G$-variety in the sense of \cite{BDP}
(recalled in Definition \ref{def reg} below). The variety $X$
contains an open orbit $G/H$ whose complement
$D$ is a strictly normal crossing divisor in $X$.
In this note we show the following vanishing result for
rational equivariant Chern classes of the bundle of logarithmic
differentials on the variety $X$:
$$
c_i^G(\Omega^1_X(\log D))=0
\quad
\text{for $i>\dim(X)-\rk(G)+\rk(H)$.}
$$

The motivation for this vanishing result 
originated in the second author's interest in a higher
rank generalization of Gieseker's proof \cite{G} of the 
Newstead-Ramanan conjecture.
The conjecture (or rather its higher rank generalization)
says that for coprime $r$ and $d$ and for $g\geq 1$ the Chern classes
of the tangent bundle of the moduli space of stable vector bundles
of rank $r$ and degree $d$ on a curve of genus $g$ vanish 
in degrees larger than $r(r-1)(g-1)$ 
(cf. \cite{Earl-Kirwan}, bottom of page 844).

In order to explain the relationship between this conjecture and
the result proven here let us first sketch Gieseker's
degeneration of moduli spaces of vector bundles.

Let $B$ a smooth curve, which serves as the base scheme
of the degeneration. Let $\X\to B$ be a proper flat family of algebraic
curves of genus $g\geq 2$
which is smooth outside a point $x\in B$.  Assume that the
fiber $X$ over $x$ is irreducible with a unique singular point $p$
which is an ordinary double point.
Let $\Xt$ be the normalization of $X$. Let $r$ be a positive integer
and let $d$ be an integer prime to $r$. 
Then there exists a variety $M(\X/B)$ proper and flat over $B$
such that
\begin{itemize}
\item
the fiber of $M(\X/B)\to B$ over any point $y\in B\setminus\{x\}$ 
is the moduli space $M(Y)$
of stable vector bundles of rank $r$ and degree $d$ on the 
curve $Y=\X_y$,
\item
the variety $M(\X/B)$ is nonsingular and its fiber $M(X)$ 
over $x$ is a normal
crossing divisor in $M(\X/B)$.
\item
Let $M(\Xt)$ be the moduli space of rank $r$ degree $d$ vector 
bundles on the curve $\Xt$. Then
there is a principal $\GL_r\times\GL_r$-bundle $P$ on $M(\Xt)$, 
and a smooth
$\GL_r\times\GL_r$-equivariant compactification
of $\GL_r=(\GL_r\times\GL_r)/\diag(\GL_r)$ which we denote
by $\KGL_r$, such that
the normalization $\Mt$ of $M(X)$ is birationally equivalent
to the locally trivial $\KGL_r$-fibration 
$$
f:M':=P\times^{\GL_r\times\GL_r}\KGL_r\to M(\Xt)
$$
associated to the principal bundle $P\to M(\Xt)$.
\end{itemize}
(cf. \cite{NS}, \cite{Se}, \cite{degeneration}, \cite{KL}).
Summarizing, we have a diagram of varieties as follows:
$$
\xymatrix{
P\times^{\GL_r\times\GL_r}\KGL_r \ar@{=}[r] &
M' \ar[d]_f \ar@{<-->}[r]^{\text{birat.}} &
\text{$\Mt$} \ar[d]^{\text{normalization}} & & \\
& \text{$M(\Xt)$} &
M(X)\ar[d] \ar@{^(->}[r] &
M(\X/B) \ar[d] &
M(Y) \ar[d] \ar@{_(->}[l] \\
& & \{x\} \ar@{^(->}[r] &
B &
\{y\} \ar@{_(->}[l]
}
$$

The conjecture of Newstead-Ramanan holds trivially for $g=1$,
since the moduli space of vector bundles of rank $r$ 
and degree $d$ on an elliptic curve $E$ is isomorphic to $E$ itself
(cf. \cite{T}). 
Gieseker's idea was to use the above diagram (in the rank two case) 
to make induction on
the genus $g$
(observe that the genus of $\Xt$ is $g-1$).

Let $D'$ be the complement of 
$
P\times^{\GL_r\times \GL_r}\GL_r
$
in $M'$ and
let $\Dt$ be the preimage of the singular locus of $M(X)$.
Then $D'$ and $\Dt$ are normal crossing divisors in $M'$ and
$\Mt$ respectively and they are proper transforms of each other by
the birational correspondence 
$
M' \leftrightarrow \Mt
$.

The induction step consists in proving the following three
implications:
\begin{eqnarray*}
& &
\text{Induction hypothesis} \\
&\overset{(1)}{=}>&
\text{Vanishing of Chern classes of $\Omega^1_{M'}(\log D')$
      in degrees larger than $r(r-1)(g-1)$} \\
&\overset{(2)}{=}>&
\text{Vanishing of Chern classes of $\Omega^1_{\Mt}(\log \Dt)$
      in degrees larger than $r(r-1)(g-1)$} \\
&\overset{(3)}{=}>&
\text{Vanishing of Chern classes of $\Omega^1_Y$
      in degrees larger than $r(r-1)(g-1)$} \\
\end{eqnarray*}

The third implication is not difficult and follows from
results already proven by Gieseker \cite{G}.
The most difficult part is the second implication. It
has been proven up to now only in the case
$r=2$ by Gieseker \cite{G}.
In the case $r=3$
Kiem and Li \cite{KL} were only able to prove the second
implication with the slightly weaker bound $6g-5$ instead of $6(g-1)$.
In both papers \cite{G} and \cite{KL} the step $(2)$ is carried through 
by means of a detailed study of flips connecting
the varieties $M'$ and $\Mt$.

It is the first implication, where the result of this note comes
to a bearing. Namely, by Example \ref{example} the
$\GL_r\times\GL_r$-equivariant embedding $\GL_r\subset\KGL_r$
is regular, and applying Corollary \ref{cor2} to this case it
follows that
$c_i(\Omega^1_{M'/M(\Xt)}(\log D'))=0$
for $i>r(r-1)$.
On the other hand, the induction hypothesis implies 
$c_i(\Omega^1_{M(\Xt)})=0$ for $i>r(r-1)(g-2)$. 
Thus the vanishing of the Chern classes 
of $\Omega^1_{M'}(\log D')$ in degrees larger than $r(r-1)(g-1)$
follows from the exact sequence
$$
0\to f^*\Omega^1_{M(\Xt)} \to \Omega^1_{M'}(\log D')
 \to \Omega^1_{M'/M(\Xt)}(\log D') \to 0
\quad.
$$

It should be noted that in the papers \cite{G} and \cite{KL}
the implication $(1)$ is proved (in the rank two and three case)
in a rather {\em ad hoc} way by
using the fact that $f:M'\to M(\Xt)$ can be represented as
a succession of blowing ups of a projective bundle.
In principle this could be done also in the higher rank case,
but as can be seen from the rank three case treated in
\cite{KL}, the computations become very involved 
with increasing rank.

\vspace{2mm}
Our Corollary \ref{cor3} establishes the vanishing in the given range
of the usual (non-equivariant) Chern classes of the bundle
of logarithmic differential forms on a regular group
embedding $G\injto X$.
This statement is proven by different methods
in a recent paper of Valentina Kiritchenko (cf. \cite{Ki} Lemma 3.6
and Proposition 4.4) where also the non-equivariant case of our
Theorem \ref{thm} is mentioned (cf. \cite{Ki} \S5).

\section{Definitions and statement of the Theorem}

\begin{definition}
Let $G$ be a topological group and let $X$ 
be a topological $G$-space.
Let $F\to X$ be a $G$-linearized complex topological 
vector bundle.
Let $E_G\to B_G$ be the universal $G$-bundle over the 
classifying space of $G$.
The {\em $G$-equivariant rational Chern class} 
$$
c_i^G(F)\in H^{2i}_G(X):=H^{2i}(E_G\times^G X,\Q)
$$
of the bundle $F$ is
by definition the Chern class of the vector bundle 
$$
E_G\times^GF\to E_G\times^GX
\quad,
$$
where $E_G\times^G F$ and $E_G\times^G X$ are
the quotients of $E_G\times F$ and $E_G\times X$ respectively
by the diagonal action of $G$.

By construction, the pull-back of $c_i^G(F)$ to $X$ (regarded as a
fiber of the map $E_G\times^G X \to B_G$) is the usual Chern class
$c_i(F)$.
\end{definition}

\begin{definition}
\label{def reg}
Let $G$ be a connected affine algebraic group and let $X$ be 
an algebraic variety on which $G$ acts with an open dense
orbit $\Omega$. The $G$-variety $X$ is called {\em regular},
if it satisfies the following conditions (see \cite{BDP}):
\begin{enumerate}
\item
The closure of every $G$-orbit is smooth.
\item
Any orbit closure $Y\neq X$ is the transversal intersection
of the orbit closures of codimension one containing $Y$.
\item
The isotropy group of any point $x\in X$ has a dense orbit in the
normal space to the orbit $G\cdot x$ in $X$.
\end{enumerate}
\end{definition}

\begin{example}
\label{example}
Let $r\geq 1$ and let $G=\GL_r\times\GL_r$.
We claim that the compactification $X:=\KGL_r$ of $\GL_r$ 
defined in \cite{kgl} 
is regular in the above sense if considered as a $G$-variety.
Indeed, properties (1) and (2) follow directly from
loc. cit. \S9. 

It remains to show property (3).
For this let $S\subseteq \GL_r$ be the maximal torus of diagonal
matrices and let $S\subset\bS_0$ be the smooth torus embedding
which in loc. cit. \S4 
we denoted by $T\subset\tilde{T}$ 
and which is an $S$-invariant
open subset of the closure $\bS$ in $X$.
From the description of the orbits in loc. cit. \S9 and the
results in loc. cit. \S4 it follows immediately that the
intersection with $\bS_0$ of any $G$-orbit of a certain codimension
in $X$ is a unique $S$-orbit of the same codimension in $\bS_0$ (or
$\bS$).

Therefore, to check property (3) it suffices to consider
points $x\in \bS_0$. Furthermore, it follows that the canonical map
$$
(N_{S\cdot x/\bS_0})_x\to (N_{G\cdot x/X})_x
$$
of normal spaces at $x$ is an $S_x$-equivariant
isomorphism.
Being a smooth torus embedding, the $S$-variety $\bS_0$
is regular in the above sense. Indeed, to see this one
is immediately reduced to the case $\Gm^n\injto\Aa^n$
where the assertion is clear.
In particular, the isotropy subgroup $S_x$
acts with a dense orbit on $(N_{S\cdot x/\bS_0})_x$. Thus 
also $G_x$ acts with a dense orbit on $(N_{G\cdot x/X})_x$.
\end{example}

Let $X$ be a regular $G$-variety. Then the complement $D\subset X$
of the open orbit in $X$ is a strict normal crossing divisor.
Recall that the sheaf of logarithmic differential forms 
$\Omega^1_X(\log D)$ is the locally free
subsheaf of $\Omega^1_X\tensor_{\Oo_X}\C(X)$
generated in a neighbourhood of a point $x\in X$ by the differentials
$$
\frac{df_1}{f_1},\dots,\frac{df_m}{f_m},df_{m+1},\dots,df_n
\quad,
$$
where $f_1,\dots,f_n$ are local coordinates at the point $x$ 
such that the divisor $D$ is given by the equation 
$f_1f_2\dots f_m=0$. As can be easily seen, this definition
is independent of the choice of the coordinate system.

The dual of $\Omega^1_X(\log D)$ is the subsheaf $T_X(-\log D)$
of the tangent sheaf $T_X$ locally at $x$ generated by the 
vector fields 
$$
f_1\frac{\partial}{\partial f_1},\dots,f_m\frac{\partial}{\partial f_m},
\frac{\partial}{\partial f_{m+1}},\dots,\frac{\partial}{\partial f_n}
\quad.
$$
In more geometric terms, it is the subsheaf of the tangent sheaf
whose sections consist of vector fields which 
are tangent to all the components of $D$.
Since the components of $D$ are $G$-invariant, it is clear from this
description that $T_X(-\log D)$ inherits a natural $G$-linearization
from the one on $T_X$. Thus $\Omega^1_X(\log D)$ is $G$-linearized
as well.

The main result of this note is the following:

\begin{theorem}
\label{thm}
Let $G$ be a connected affine algebraic group over the field of
complex numbers and let $X$ be a complex regular $G$-variety.
Let $H\subseteq G$ be the isotropy group of a point in the open orbit
$\Omega\subseteq X$ and let $D:=X\setminus\Omega$ be the boundary divisor.
Then the $G$-equivariant rational Chern classes of the bundle
$\Omega_X^1(\log D)$ of logarithmic differentials vanish in degrees
larger than 
$$
\dim(X)-\rk(G)+\rk(H)
\quad .
$$
\end{theorem}

The proof will be given in the next section.
We note that the corresponding statement over an algebraically
closed field of arbitrary characteristic holds true as well
with the same proof 
if one replaces rational equivariant cohomology with the equivariant
Chow ring with rational coefficients. 

The following is an immediate consequence:

\begin{corollary}
\label{cor1}
Let $G$ be a complex connected affine algebraic group and let
$G=(G\times G)/\diag(G)\injto X$ be a regular 
$G\times G$-equivariant embedding
with boundary divisor $D$. 
Then we have $c_i^{G\times G}(\Omega^1_X(\log D))=0$ for
$i>\dim(G)-\rk(G)$.
\end{corollary}

For the application mentioned in the introduction, we state the following
\begin{corollary}
\label{cor2}
Let $G$, $X$ and $D$ be as in \ref{cor1}.
Let $M$ be a complex variety and let $P$ be a principal
$G\times G$-bundle over $M$. Let 
$$
\X:=P\times^{G\times G}X 
\quad\text{and}\quad
\D:=P\times^{G\times G}D \subset\X
$$
be the associated locally trivial $X$-fibration and $D$-fibration
over $M$ respectively. 
Then the Chern classes of $\Omega^1_{\X/M}(\log \D)$ vanish in degrees
larger than $\dim(G)-\rk(G)$.
\end{corollary}

\begin{proof}
Let $M\to B_{G\times G}$ be a classifying map for the bundle $P$.
Then we have cartesian diagrams
$$
\xymatrix{
\text{$\X$} \ar[d]\ar[r]^(.3)f & \text{$E_{G\times G}\times^{G\times G}X$} \ar[d] \\
M \ar[r] &  \text{$B_{G\times G}$}
}
\qquad\qquad
\xymatrix{
\text{$\D$} \ar[d]\ar[r] & \text{$E_{G\times G}\times^{G\times G}D$} \ar[d] \\
M \ar[r] &  \text{$B_{G\times G}$}
}
$$
and $\Omega^1_{\X/M}(\log \D)$ is the pull back via
$f$
of the vector bundle 
$$
E_{G\times G}\times^{G\times G}\Omega^1_{X}(\log D)
\quad.
$$
Therefore we have
$$
c_i(\Omega^1_{\X/M}(\log \D))=
f^*c_i(E_{G\times G}\times^{G\times G}\Omega^1_{X}(\log D))=
f^*c_i^{G\times G}(\Omega^1_{X}(\log D))
$$
and the result follows from Corollary \ref{cor1}.
\end{proof}

\begin{corollary}
\label{cor3}
With the notation of Corollary \ref{cor1}, the usual Chern classes of
$\Omega^1_X(\log D)$ vanish in degrees $> \dim(G)-\rk(G)$.
\end{corollary}

\begin{proof}
This is a special case of Corollary \ref{cor2}: Take $M$ to be a point.
\end{proof}

\section{Proof of the theorem}

The first step in the proof 
of the theorem is the case of empty boundary divisor:

\begin{proposition}
\label{G/H}
Theorem \ref{thm} holds in the case when $X=G/H$.
\end{proposition}

\begin{proof}
We have 
to show the vanishing of the $G$-equivariant Chern 
classes of the cotangent bundle
$\Omega_X^1$ in the given range.
Let $x\in X$ be the point represented by the unit element in $G$.
The fiber $V$ of $\Omega_X^1$ at $x$ is an $H$-module. In fact we have
$$
V=(\Lie(G)/\Lie(H))^{\dual}
\quad,
$$
where the action of $H$ on the right hand side is 
induced by the adjoint action of $H$ on $\Lie(G)$.

Since $H$ operates freely on the contractible space $E_G$,
we may take $E_H=E_G$ and $B_H=E_G/H=E_G\times^GX$.
Since $\Omega^1_X=G\times^HV$, it follows that the vector bundle
$E_G\times^G\Omega^1_X\to E_G\times^GX$ can be identified with
the vector bundle 
$$
\V:=E_H\times^HV\to B_H
\quad.
$$
Thus we have to show the vanishing 
of the Chern classes of $\V$ 
in degrees larger than
$$
\dim(X)-\rk(G)+\rk(H)
\quad .
$$

Let $S_H$ be a maximal torus of $H$ and let $X^*(S_H)$ be its character
group. By the results of \cite{H} Chapter III, \S1, the map
$B_{S_H}\to B_H$ induces an injection
$$
H^*(B_H,\Q)\injto H^*(B_{S_H},\Q)\isomto\Sym_{\Q}(X^*(S_H)\tensor\Q)
\quad.
$$
Now if we denote by $\alpha_1,\dots,\alpha_{\dim X}\in X^*(S_H)$ 
the weights of $V$ considered as an $S_H$-module,
then the total Chern class of $\V$ is mapped to the product
$$
\prod_{i=1}^{\dim X}(1+\alpha_i)\in \Sym_{\Q}(X^*(S_H)\tensor\Q)
$$

But $S_H$ acts trivially on the subspace
$$
\Lie(S_G)/\Lie(S_H) \subseteq \Lie(G)/\Lie(H)=V^{\dual}
\quad,
$$
where $S_G$ denotes a maximal torus of $G$ containing $S_H$.
Therefore at least as many as $\dim(\Lie(S_G)/\Lie(S_H))$ of the 
$\alpha_i$ vanish and
the total Chern class of $\V$ is of degree
at most 
$$
\dim(X) -\dim(\Lie(S_G)/\Lie(S_H))=\dim(X)-\rk(G)+\rk(H)
\quad.
$$
\end{proof}

\begin{proposition}
\label{step}
Let $G$ be a connected affine algebraic group and let $X$ be
a regular $G$-variety.
Let $x,y\in X$ be two points and let
$G\cdot x$ and $G\cdot y$ be the corresponding $G$-orbits in $X$.
Assume
that the orbit closure $\overline{G\cdot y}$ is contained 
in $\overline{G\cdot x}$ as a divisor.
Let $G_x$ and $G_y$ be the isotropy subgroups of the
points $x$ and $y$.
Then
$$
\rk(G_x)\geq \rk(G_y)-1
\quad.
$$ 
\end{proposition}

\begin{proof}
Since $\overline{G\cdot x}$ is again a regular $G$-variety,
we may assume 
without loss of generality that $\overline{G\cdot x}=X$.
This simplifies notation a bit.

Let $S_y$ be a maximal subtorus of $G_y$.
It acts on the tangent space $T_y(X)$ and leaves
the codimension-one subspace $T_y(G\cdot y)$ invariant.
Therefore there is a one-dimensional $S_y$-invariant subspace $\ell$ in
$T_y(X)$ such that 
$$
T_y(X)=T_y(G\cdot y)\oplus \ell
\quad.
$$

Let $\chi:S_y\to \C^{\times}$ be the character
associated to the one-dimensional $S_y$-module $\ell$
and let $S$ be the connected component of the identity in $\ker(\chi)$.
We claim that there is a point $x_0\in G\cdot x$
which is fixed by $S$.

Indeed, since $X$ is normal, by a result of Sumihiro (\cite{Su})
there is an affine open neighbourhood $U\isomorph\Spec(R)$ of $y$ in $X$
which is invariant under the action of $S$.
Let $\m\subset R$ be the maximal ideal corresponding to the point 
$y\in U$.
Since $S$ is a torus, there is an $S$-equivariant section
$$
\varphi:\m/\m^2\to \m
$$
of the surjection $\m\to\m/\m^2$.
The map $\varphi$ induces an $S$-equivariant
map $\Sym(\m/\m^2)\to R$ and thus an $S$-equivariant morphism
$$
f: U=\Spec(R)\to\Spec(\Sym(\m/\m^2))=T_y(X)
\quad.
$$
By construction, $f$ induces an isomorphism from the
tangent space of $U$ at $y$ to the tangent space of $T_y(X)$ at $0$.
Therefore after possibly restricting $f$ to a smaller $S$-invariant
neighbourhood of $y$, we may assume that $f$ is \'etale.

Let $Z$ be the preimage by $f$ of the line $\ell$.
Since $\ell$ is transversal to the tangent directions
of $G\cdot y$, it follows that $Z$ is not contained
in $G\cdot y$. Let $x_0\in Z\setminus G\cdot y$.

Let us collect what we know about the point $x_0$.
The image $f(x_0)$ of $x_0$ in $T_y(X)$ is contained in $\ell$ and
is thus fixed by $S$. Since $f$ is $S$-equivariant, the group
$S$ acts on the fiber of $f$ at $f(x_0)$. But this fiber contains
only finitely many points and $S$ is connected, so $S$ must act trivially
on the fiber. In particular, $x_0$ is a fixed point of $S$ contained
in $G\cdot x$. This proves our claim.

It follows that $S\subseteq G_{x_0}$ and consequently we have
$$
\rk(G_x)=\rk(G_{x_0})\geq \dim(S)\geq \dim(S_y)-1 = \rk(G_y)-1
\quad.
$$
\end{proof}

\begin{corollary}
\label{inequality}
Let $G$ be a connected affine algebraic group and let
$X$ be a regular $G$-variety.
Let $H\subset G$ be the isotropy group of a point in the open orbit.
Then for any point $y\in X$ with isotropy group $G_y$
we have the inequality
$$
\dim(H)-\rk(H)\leq \dim(G_y) - \rk(G_y)
\quad.
$$
\end{corollary}

\begin{proof}
Let $c$ be the codimension in $X$ of the orbit $G\cdot y$.
There is a sequence of points $x_0,x_1,\dots,x_c$ in $X$ 
such that $x_0$ is contained in the open orbit and $x_c=y$, and
such that for $i=1,\dots,c$ the orbit closure $\overline{G\cdot x_i}$ 
is a divisor in
$\overline{G\cdot x_{i-1}}$. Applying Proposition 
\ref{step} $c$ times, it follows that
$$
\rk(H)\geq\rk(G_y)-c
\quad.
$$

On the other hand we have 
$$
\dim(G)-\dim(H)-c=\dim(X)-c=\dim(G\cdot y)=\dim(G)-\dim(G_y)
$$
and therefore $\dim(H)=\dim(G_y)-c$.
From this the inequality of Corollary \ref{inequality} is immediate.
\end{proof}

\begin{proposition}
\label{exact}
Let $G$ be a connected affine algebraic group and let $X$ be
a regular $G$-variety. Let $D\subset X$ be the boundary divisor.
Let $Y=G\cdot y$ be a $G$-orbit of codimension
$c$ in $X$.
Then there is an exact sequence of $G$-linearized
bundles on $Y$ as follows:
$$
0\to \Omega^1_Y \to \Omega^1_X(\log D)|_Y \to \Oo_Y^{\oplus c} \to 0     
\quad.
$$
\end{proposition}

\begin{proof}
This follows from \cite{BB}, 2.4.2.
\end{proof}

Now we can prove Theorem \ref{thm}.
For each $G$-orbit $Y\subseteq X$ we have
the restriction map $H^*_G(X)\to H^*_G(Y)$.
Thus we get a homomorphism
$$
H^*_G(X)\to\prod_{Y}H^*_G(Y)
$$
where $Y$ runs through all $G$-orbits of $X$.
By \cite{BDP} Theorem 7 this map is injective.

Therefore it suffices to show that for any
$$
i>\dim(G)-\dim(H)-\rk(G)+\rk(H)
$$
and any $G$-orbit $Y=G\cdot y$ in $X$ the
$i$-th $G$-equivariant rational 
Chern class of $\Omega^1_X(\log D)|_Y$ vanishes.
By Proposition \ref{exact} this is equivalent to the vanishing
of the $i$-th $G$-equivariant rational Chern class of $\Omega^1_Y$.

Let $G_y$ be the isotropy group of the point $y$.
From Corollary \ref{inequality} we know that
$$
i>\dim(G)-\dim(H)-\rk(G)+\rk(H)\geq
  \dim(G)-\dim(G_y)-\rk(G)+\rk(G_y)
\quad.
$$
Proposition \ref{G/H} yields now the vanishing of
$c^G_i(\Omega^1_Y)$ as required.


\begin{thebibliography}{99}
\bibitem[BB]{BB}
        F. Bien and M. Brion:
        Automorphisms and local rigidity of regular varieties.
        Compositio Mathematica 104, 1--26, 1996.
\bibitem[BDP]{BDP}
        E. Bifet, C. De Concini, C.Procesi:
        Cohomology of Regular Embeddings,
        Advances in Mathematics 82, 1--34 (1990)
\bibitem[EK]{Earl-Kirwan}
        R. Earl, F. Kirwan:
        The Pontryagin rings of moduli spaces of arbitrary rank 
        holomorphic bundles over a Riemann surface. 
        J. London Math. Soc. (2) 60 (1999), no. 3, 835--846. 
\bibitem[G]{G}
        D. Gieseker:
        A degeneration of the moduli space of stable bundles.
        J. Differential Geometry 19 (1984) 173--206.
\bibitem[H]{H}
        W. Y. Hsiang:
        Cohomological Theory of Topological Transformation Groups,
        Ergeb. der Math. 85, Springer-Verlag, New York, 1975.  
\bibitem[K1]{kgl}
        I. Kausz:
        A modular compactification of the general linear group.
        Documenta Mathematica 5 (2000) 553--594.
\bibitem[K2]{degeneration}
        I. Kausz:
        A Gieseker type degeneration of moduli stacks of vector 
        bundles on curves.
        Transactions of the AMS, to appear.
\bibitem[KL]{KL}
        Y.-H. Kiem and J. Li:
        Vanishing of the top Chern classes of the moduli of
        vector bundles.
        Preprint math.AG/0403033.
\bibitem[Ki]{Ki}
        V. Kiritchenko:
        Chern classes of compactifications of reductive groups.
        Preprint math.AG/0411331.
\bibitem[NS]{NS}
        D.S. Nagaraj and C.S. Seshadri:
        Degenerations of the moduli spaces of vector bundles on 
        curves II.
        Proc. Indian Acad. Sci. Math. Sci. 109 (1999), no 2, 165--201.
\bibitem[Se]{Se}
        C.S. Seshadri: 
        Degenerations of the moduli spaces of vector bundles on curves. 
        School on Algebraic Geometry (Trieste, 1999), 205--265, 
        ICTP Lect. Notes, 1, 
        Abdus Salam Int. Cent. Theoret. Phys., Trieste, 2000. 
\bibitem[Su]{Su}
        H. Sumihiro:
        Equivariant completion. 
        J. Math. Kyoto Univ. 14 (1974), 1--28.
\bibitem[T]{T}
        L.W. Tu:
        Semistable Bundles over an Elliptic Curve.
        Advances in Mathematics 98, 1--26 (1993).
        
\end{thebibliography}
\end{document}